\documentclass[11pt,twoside]{amsart}

\usepackage{amssymb}
\usepackage{latexsym}

\newcommand{\Id}{\mathrm{Id}}      % Copy editor: This is *NOT* a
\newcommand{\dopu}{{:}\allowbreak\ }
\newcommand{\eps}{\varepsilon}
\newcommand{\cal}{\mathcal}

\newcommand{\loglike}[1]{\mathop{\rm #1}\nolimits}

\newcommand{\coq}{\loglike{\overline{co}}}

\newcommand{\dist}{\loglike{dist}}

\newcommand{\re}{\loglike{Re}}
\newcommand{\im}{\loglike{Im}}

\newcommand{\N}{{\mathbb N}}

\newcommand{\C}{{\mathbb C}}
\newcommand{\T}{{\mathbb T}}
\newcommand{\Z}{{\mathbb Z}}

\theoremstyle{plain}
\newtheorem{thm}{Theorem}[section]
\newtheorem{theo}[thm]{Theorem}

\newtheorem{prop}[thm]{Proposition}
\newtheorem{cor}[thm]{Corollary}
\newtheorem{lemma}[thm]{Lemma}

\theoremstyle{definition}
\newtheorem{definition}[thm]{Definition}

\theoremstyle{remark}

\numberwithin{equation}{section}

\newcommand{\rest}[2]{#1\raisebox{-0.3ex}{\mbox{$\mid_{#2}$}}}
\def\DP{Daugavet property}

\newcommand{\begsta}{\begin{statements}}
\newcommand{\begaeq}{\begin{aequivalenz}}
\def\endsta{\end{statements}}
\def\endaeq{\end{aequivalenz}}
\newcommand{\bea}{\begin{eqnarray*}}
\newcommand{\eea}{\end{eqnarray*}}

\newcounter{abc}   % Counter f\"ur statements-environment wird deklariert
\newcounter{iiiii} % Counter f\"ur aequivalenz-environment wird deklariert

\newenvironment{aequivalenz}
{\setcounter{iiiii}{0}
\begin{list}%
{{\rm (\roman{iiiii})}}%  Falls die items nicht angegeben sind: i)u.s.w.
{\usecounter{iiiii}
\parsep=0pt plus 1pt
\topsep=1pt plus 2pt minus 1pt
\itemsep=1pt plus 2pt minus 1pt
\leftmargin=3\baselineskip
\labelsep=.6\baselineskip
\labelwidth=2.4\baselineskip
\rightmargin 0pt}%
}%
{\end{list}}

\newenvironment{statements}%
{\setcounter{abc}{0}
\begin{list}%
{{\rm (\alph{abc})}}%  Falls die items nicht angegeben sind: (a) u.s.w.
{\usecounter{abc}
\parsep=0pt plus 1pt
\topsep=1pt plus 2pt minus 1pt
\itemsep=1pt plus 2pt minus 1pt
\leftmargin=3\baselineskip
\labelsep=.6\baselineskip
\labelwidth=2.4\baselineskip
\rightmargin 0pt}%
}%
{\end{list}}

\begin{document}

\title{ Remarks on rich subspaces of Banach spaces }

\author{ Vladimir Kadets, Nigel Kalton and Dirk Werner }

\address{Faculty of Mechanics and Mathematics, Kharkov National University,
\qquad {}\linebreak pl.~Svobody~4,  61077~Kharkov, Ukraine}
%\email{anna.m.vishnyakova@univer.kharkov.ua}
\email{vova1kadets@yahoo.com}

\curraddr{Department of Mathematics, Freie Universit\"at Berlin,
\mbox{Arnimallee~2--6}, D-14\,195~Berlin, Germany}
\email{kadets@math.fu-berlin.de}

\address{Department of Mathematics, University of Missouri,
Columbia MO 65211}
\email{nigel@math.missouri.edu}

\address{Department of Mathematics, Freie Universit\"at Berlin,
Arnimallee~2--6, \qquad {}\linebreak D-14\,195~Berlin, Germany}
\email{werner@math.fu-berlin.de}

%\date{\today}

\thanks{The work of the first-named author was supported by a
fellowship from the Alexander von Humboldt Foundation. 
The second-named author was supported by NSF grant DMS-9870027.}



\begin{abstract}
We investigate rich subspaces of $L_1$ and deduce an
interpolation property of Sidon sets.
 We also present examples of rich
separable subspaces of nonseparable Banach spaces and we study the
Daugavet property of tensor products.
\end{abstract}

\maketitle

%%%%%%%%%%%%%%%%%%%%%%%%%%%%%%%%%%%%%%%%%%%%
\thispagestyle{empty}

\begin{center} \em
Dedicated to Professor Aleksander Pe{\l}czy\'nski \\
on the occasion of his 70th birthday \\[1cm]
\end{center}

\section{Introduction}

In this paper we present some results concerning the notion of a rich
subspace of a Banach space as introduced in \cite{KadSW2}. In that
paper
(see also \cite{dirk-irbull}),
an operator $T\dopu X \to Y$ is called \textit{narrow} if for every $x, y  \in
S(X)$ (the unit sphere of $X$),
$\eps > 0$  and every slice $S $
of the unit ball $B(X)$ of $X$ containing $y$ there is an element $v \in S$
such that  $\|x+v\|>2-\eps$ and
$\|T(y-v)\|< \eps$, and a subspace $Z$ of $X$ is called \textit{rich} if
the quotient map $q\dopu X\to X/Z$ is narrow. 
We recall that a slice of the unit ball is a nonvoid set of the form
$S= \{x\in B(X)\dopu \re x^*(x) > \alpha\}$
for some functional $x^*\in X^*$. 
Thus, $Z$ is a rich subspace if for every $x, y  \in
S(X)$,
$\eps > 0$  and every slice $S $
of $B(X)$ containing $y$ there is some $z\in X$ at distance ${\le
\eps}$ from $Z$ such that $y+z\in S$ and $\|x+y+z\|> 2-\eps$. 
Actually, we are not giving the original definition 
of a narrow operator but the equivalent
reformulation from \cite[Prop.~3.11]{KadSW2}.

These ideas build on
previous work in \cite{PliPop} and  \cite{KadPop}; however we point
out that the above definition of richness is unrelated to Bourgain's
in \cite{Bour-84b}. 
Narrow operators were used in \cite{BKSSW} and \cite{KadPop}
to extend Pe{\l}czy\'nski's classical result that neither $C[0,1]$
nor $L_{1}[0,1]$ embed into spaces having unconditional bases.

The investigation of narrow operators is closely connected with the
\DP\ of a Banach space. A Banach space $X$ has the \textit{\DP} whenever
$\|\Id + T\| = 1 + \|T\|$ for every rank-1 operator $T\dopu X\to
X$; prime examples are $C(K)$ when $K$ is perfect (i.e., has no
isolated points), $L_{1}(\mu)$ and $L_{\infty}(\mu)$ when $\mu$ is
nonatomic, the
disc algebra, and spaces like $L_{1}[0,1]/V$ when $V$ is reflexive. 
For future reference we mention the following characterisation of
the \DP\ \cite{KadSSW}:

\begin{lemma} \label{l:MAIN}
The following assertions are equivalent:
\begin{aequivalenz}
\item
     $X$ has the Daugavet property.
\item
     For every $x\in S(X)$, $\eps>0$ and every slice $S$
of $B(X)$ there exists some $v\in S$ such that $\|x+v\|>2-\eps$.
\item
For all $x\in S(X)$ and $\eps>0$, $B(X) = \coq \{v\in B(X)\dopu
\|x+v\|>2-\eps\}$.
\end{aequivalenz}
\end{lemma}

Therefore, $X$ has the \DP\ if and only if $0$ is a narrow operator 
on $X$ or equivalently  if
and only if there exists at least one narrow operator on~$X$. 
It is proved in \cite{KadSW2} that then every weakly compact operator 
on $X$ with values in some Banach space $Y$ (indeed,
every strong Radon-Nikod\'ym operator) and every operator not fixing a
copy of $\ell_1$ is narrow (and hence satisfies $\|\Id + T\| = 1 +
\|T\|$ when it maps $X$ into $X$). Consequently, a subspace $Z$ of a
space with the \DP\ is rich if $X/Z$ or $(X/Z)^*$ has the RNP.

Also, $X$ has the \DP\ if and only if $X$ is a rich subspace in itself
or equivalently if $X$ contains at least one rich subspace.

The general idea of
these notions is that a narrow operator is sort of small and hence a
rich subspace is large. In Section~\ref{sec2} of this paper we study
rich subspaces of $L_1$. With reference to a quantity that is
reminiscent of the Dixmier characteristic we show that a rich subspace
is indeed large: a subspace with a bigger ``characteristic'' coincides
with $L_1$. As an application we present an interpolation property of
Sidon sets. 
We remark that the counterpart notion of a small subspace of $L_1$ has
been defined and investigated in \cite{GoKaLi2}. 

These results  notwithstanding, Section~\ref{sec3}  gives examples of
rich subspaces that appear to be small, namely there are   examples of
nonseparable spaces and separable rich subspaces. 

In Section~\ref{sec4} we study hereditary properties for the \DP\ in
tensor products. Although there are positive results for rich
subspaces of $C(K)$, we present counterexamples in the general case.

%%%%%%%%%%%%%%%%%%%%%%%%%%%%%%%%%%%%%%%%%%%%%%%%%%%%%%%%%%%%

\section{Rich subspaces of $L_{1}$} 
\label{sec2}

Let $X\subset L_1 = L_1(\Omega,\Sigma,\lambda)$ be a closed subspace
where $\lambda$ is a probability measure.
We define $C_X$ to
be the closure of $B(X)$ in $L_1$ with respect to the $L_0$-topology,
the topology of convergence in measure. Note that for $f\in C_X$
there is a sequence $(f_n)$ in $B(X)$ converging to $f$ pointwise almost
everywhere and almost uniformly. In this section, the symbol
${\|f\|}$ refers to the $L_1$-norm of a function.

In \cite[Th.~6.1]{KadSW2} narrow operators  on the 
space $L_1$ were characterised as follows.

\begin{theo} \label{6.1}
An operator $T\dopu L_1 \to Y$ is narrow if and only if for every
measurable set $A$ and every $\delta,\eps >0$ there is a real-valued
$L_1$-function $f$ supported on $A$ such that $\int f=0$, $f\le 1$,
the set $\{f=1\}$ of those $t\in\Omega$ for which $f(t)=1$
has measure $\lambda ( \{f=1\} ) > \lambda(A)-\eps$ and $\|Tf\|\le \delta$.
In particular, a subspace $X\subset L_{1}$ is rich if and only if
for every
measurable set $A$ and every $\delta,\eps >0$ there is a real-valued
$L_1$-function $f$ supported on $A$ such that $\int f=0$, $f\le 1$,
$\lambda(\{f=1\}) > \lambda(A)-\eps$ and the distance from $f$
to $X$ is ${\le \delta}$.
\end{theo}

Actually, in \cite{KadSW2} only the case of real $L_1$-spaces was
considered, but the proof extends  to the complex case. Indeed,
instead of the function $v$ that is constructed in the first part of
the proof of \cite[Th.~6.1]{KadSW2} one uses its real part and employs
the fact that for real-valued $L_1$-functions $v_1$ and $v_2$
satisfying
$$
1-\delta < \int_{\Omega} |v_1| \, d\lambda \le
\int_\Omega (v_1^2 + v_2^2)^{1/2} \,d\lambda \le 1
$$
we have $\|v_2\|\le  \sqrt{2\delta} $.

%We will use this theorem repeatedly.

\begin{prop} \label{prop2}
If $X$ is rich, then $\frac12 B({L_1}) \subset C_X$.
\end{prop}

\begin{proof}
Since $C_X$ is $L_1$-closed, it is enough to show that $f_A := 
\chi_A/\lambda(A) \in 2C_X$ for every measurable set $A$. By
Theorem~\ref{6.1} there is, given $\eps>0$, a 
real-valued function $g_\eps$
supported on $A$ with $g_\eps \le 1$ and $\int g_\eps =0$ such that
 $\{ g_\eps <1 \} $ has measure $\le\eps$ and the distance of $g_\eps$
to $X$ is $\le\eps$. Clearly $g_\eps/\lambda(A) \to f_A$ in measure
as $\eps\to0$ and
$$
\|g_\eps\| = \|g_\eps^+\| + \|g_\eps^{-}\| = 2 \|g_\eps^+\| \le 2\lambda(A).
$$
Therefore, there is a sequence $(f_n)$ in $X$ of norm ${\le2}$
converging to $f_A$ in measure.
\end{proof}

\begin{prop} \label{prop3}
If $\frac12 B({L_1}) \subset C_Y$ for all $1$-codimensional subspaces
$Y$ of $X$, then $X $ is rich.
\end{prop}

\begin{proof}
Again by Theorem~\ref{6.1}, we have to produce functions 
$g_\eps$ as above on any given measurable set~$A$.
Therefore, we let $Y= \{ f\in X\dopu \int_A f =0 \} $.
By assumption, there is a sequence $(f_n)$ in $Y$ such that 
$\|f_n\|\le2 \lambda(A)$ and $f_n\to \chi_A$ in measure.

We shall argue that $\|{\im f_n}\|\to 0$. 
Let $\eta >0$. If $n$ is large enough, the set $B_n:= \{|f_n-\chi_A| \ge \eta \} $
has measure ${\le\eta}$. For those~$n$,
$$
0= \int_{A} {\re f_n} =
\int_{A\setminus B_n} {\re f_n} + \int_{A\cap B_n} {\re f_n} 
$$
implies that
$$
\int_{A\cap B_n} |{\re f_n}| \ge 
\biggl| \int_{A\cap B_n} {\re f_n} \biggr|  =
\biggl| \int_{A\setminus B_n} {\re f_n} \biggr|  \ge
\lambda(A\setminus B_n) (1-\eta)
$$
and
$$
\|{ \rest{\re f_n}{A} } \| \ge
\lambda( A\setminus B_{n}) (1-\eta) + 
 \int_{A\cap B_n} |{\re f_n}| \ge
2(\lambda(A)-\eta) (1-\eta).
$$
Hence,
$$
2(\lambda(A ) - \eta) (1-\eta) \le
\|{ \rest{\re f_n}{A} } \| \le
\|{ \rest{ f_n}{A} } \| \le
\| f_n \| \le
2 \lambda(A),
$$
and it follows for one thing that $\|{ \rest{\im f_n}{A} } \|$ is small provided $\eta$ is 
small enough (cf.\ the remarks after Theorem~\ref{6.1}) and moreover that
$$
\|{ \rest{ f_n}{[0,1]\setminus A} } \| \le
2\eta + 2\eta \lambda(A).
$$
Consequently, $\|{ {\im f_n} } \| \to 0$
as $n\to \infty$.

Now let $\delta= \eps/9$ and choose $n$ so  large that the set $B:= \{ |{\re f_n-\chi_A|}
\ge\delta \} $ has measure ${\le\delta}$ and $\|{ \im f_n } \| \le \delta$. 
Then there exists a 
real-valued function 
$h$ such that $h=0$ on $[0,1]\setminus (A\cup B)$, $h=1$ on $A \setminus B$,
$\int_A h =0$ and $\|h- \re f_n\|\le 2\delta$. Now
\begin{align}
\| \rest{h}{A} \| &=
2\| \rest{h^+}{A} \| \ge 2(\lambda(A)-\delta) \nonumber \\
\|h\| &\le 
\|{\re f_n}\| + 2\delta \le 2(\lambda(A)+\delta), \nonumber
\end{align}
so
$$
\| \rest{h}{[0,1]\setminus A} \| \le 4 \delta.
$$
Furthermore,
\begin{align}
\| \rest{h^+}{A} \| &=
\| \rest{h^+}{A\cap B} \| + \| \rest{h^+}{A\setminus B} \|
\ge \| \rest{h^+}{A\cap B} \| + \lambda(A) - \delta, 
\nonumber \\
2 \| \rest{h^+}{A} \| &=
\| \rest{h}{A} \| \le 2(\lambda(A) + \delta), \nonumber
\end{align}
so
$$
\| \rest{h^+}{A\cap B} \| \le 2\delta,
$$
and it follows that there is a function $g=g_\eps$ such that $g=0$ on
$[0,1]\setminus A$, $g=1$ on $A\setminus B$, $\int g = 0$, $g\le1$
and $\|g-h\|\le 4\delta$. Then
$$
\dist(g,X)\le \|g-f_n\| \le \|g-h\|+ \|h- \re f_n\| + \|{\im f_n}\| 
\le 9\delta = \eps, 
$$
as requested.
\end{proof}

Since a $1$-codimensional subspace of a rich subspace is rich
\cite[Th.~5.12]{KadSSW},
Proposition~\ref{prop2} shows that Proposition~\ref{prop3} can
actually be formulated as an equivalence. This is not so for
Proposition~\ref{prop2}: the space constructed in Theorem~6.3
of \cite{KadSW2} is not rich, yet it satisfies $\frac12 B({L_1}) \subset C_X$.

We sum  this up in a theorem.

\begin{theo} \label{theo4}
$X$ is a rich subspace of $L_{1}$ if and only if
$\frac12 B({L_1}) \subset C_Y$ for all $1$-codimensional subspaces
$Y$ of $X$.
\end{theo}

The next proposition shows that the factor $\frac12$ is optimal.

\begin{prop} \label{prop4}
If, for some $r>\frac12$, $rB({L_1})\subset C_X$, then $X=L_1$.
\end{prop}

\begin{proof}
Suppose $h\in L_\infty$, $\|h\|_\infty =1$, and let $Y= \{ f\in L_1\dopu
\int fh =0 \}$. Assume that $B({L_1})\subset s C_Y$; we shall argue that
$s\ge2$. This will prove the proposition since every proper closed
subspace is contained in a closed hyperplane.

Assume without loss of generality that $h$ takes the (essential)
value~$1$. 
Let $\eps>0$, and put $A=\{ |h-1| <\eps/2 \} $; then $A$ has positive
measure. There is a sequence $(f_n)$
converging to $\chi_A$ in measure such that $\|f_n\|\le s\,\lambda(A)$ and
$\int f_n h =0$ for all~$n$. Since $f_nh\to \chi_A h$ in measure as well,
there is, if $n$ is a sufficiently large index,
 a subset $A_n\subset A$ of measure ${\ge (1-\eps)\lambda(A)}$
such that $|f_nh - 1| <\eps$ on~$A_{n}$.
For such an~$n$, 
\bea
\biggl| \int_{A_{n}} f_{n}h  \biggr|
&=&
\biggl| \lambda(A_{n}) - \int_{A_{n}} (1-f_{n}h) \biggr| \\
&\ge&
\lambda(A_{n}) - \int_{A_{n}} |1-f_{n}h| 
\ge
(1-\eps) \lambda(A_{n})  ,
\eea
and therefore
$$
\int_{A_{n}} |f_{n}h  | \ge (1-\eps) \lambda(A_{n})
$$
and, if $B_{n}$ denotes the complement of $A_{n}$,
$$
\int_{B_{n}} |f_{n}h  | \ge 
\biggl| \int_{B_{n}} f_{n}h  \biggr| =
\biggl| \int_{A_{n}} f_{n}h  \biggr| \ge 
(1-\eps) \lambda(A_{n})
$$
so that
$$
s\, \lambda(A) \ge \|f_n\| \ge \|f_nh\| \ge 
 2(1-\eps)^2 \lambda(A).
$$
Since $\eps>0$ was arbitrary, we conclude that $s\ge2$.
\end{proof}

Thus, the rich subspaces appear to be the next best thing in terms
of size of a subspace after $L_1$ itself. At the other end of the 
spectrum are the nicely placed subspaces, defined by the condition
that $B(X)$ is $L_0$-closed. Recall that $X$ is nicely placed if
and only if $X$ is an $L$-summand in its bidual, i.e., 
$X^{**}= X \oplus_1 X_s$ ($\ell_1$-direct sum) for some closed
subspace $X_s$ of $X^{**}$ \cite[Th.~IV.3.5]{HWW}.

We now look at the translation invariant case, and 
we consider $L_1(\T)$ (or $L_1(G)$ for a compact %metrisable
abelian group). As usual,
for $\Lambda\subset \Z$ the space $L_{1,\Lambda}$ 
consists of those $L_1$-functions 
whose Fourier coefficients vanish off~$\Lambda$.

\begin{prop} \label{prop5}
Let $\Lambda\subset \Z$ and suppose that $L_{1,\Lambda}$ is rich in
$L_{1}$. Then for every measure $\mu$ on $\T$ and every $\eps>0$ 
there is a measure $\nu$ with $\|\nu\| \le \|\mu\| + \eps$ 
and $\widehat{\nu}(\gamma) = \widehat{\mu}(\gamma)$ 
for all $\gamma\notin \Lambda$ that is $\eps$-almost singular in the
sense that there is a set $S$ with $\lambda(S)\le \eps$ and
$|\nu|(\T \setminus S) \le \eps$.
\end{prop}

\begin{proof}
Let $\mu = f\lambda + \mu_{s}$ be the Lebesgue decomposition of
$\mu$, and let $\delta>0$. By Proposition~\ref{prop2}
there is a function $g\in L_{1,\Lambda}$ such that
$\|g\|\le 2\|f\|$ and $A:=\{|f-g|>\delta\}$ has measure 
${<\delta}$. Let $B:=\{|f-g|\le\delta\}$. Then
$$
\|g\chi_A\| \le 2\|f\| - \|g\chi_B\| \le
2\|f\| - \|f\chi_B\| + \delta =
\|f\| + \|f\chi_A\| + \delta.
$$
Therefore we have for $\nu := \mu - g\lambda$
\bea
\|\nu\| &=&
\|(f-g)\lambda + \mu_s\| \\
&\le&
\|f\chi_A\| + \|g\chi_A\| + \|(f-g)\chi_B\| + \|\mu_s\| \\
&\le&
2\|f\chi_A\| + 2\delta + \|\mu\|,
\eea
and hence $\|\nu\|\le \|\mu\| + \eps$ if $\delta$ is sufficiently
small.

Clearly $\widehat{\nu} = \widehat{\mu}$ on the complement of
$\Lambda$, and if $N$ is a null set supporting $\mu_{s}$, then $S:= A
\cup N$ has the required properties if $\delta\le\eps$.
\end{proof}

We apply these ideas to Sidon sets, i.e., sets $\Lambda'\subset \Z$
such that all functions in $C_{\Lambda'}$ have absolutely sup-norm
convergent Fourier series. (See \cite{KalPel} for recent results on this
notion.)
If $\Lambda$ is the complement of a Sidon set, then $L_1/L_{1,\Lambda}$
is isomorphic to $c_0$ or finite-dimensional \cite[p.~121]{Rud-F}. 
Hence $L_{1,\Lambda}$ is rich by \cite[Prop.~5.3]{KadSW2},
and Proposition~\ref{prop5}  applies.
Thus, the following corollary holds.

\begin{cor} \label{cor6}
If $\Lambda'\subset \Z$ is a Sidon set and
$\mu$ is a measure on $\T$, then for every $\eps>0$
there is an $\eps$-almost singular measure $\nu$
with $\|\nu\| \le \|\mu\|+\eps$ and $\widehat{\nu}(\gamma) = 
\widehat{\mu}(\gamma)$ for all $\gamma\in \Lambda'$.
\end{cor}

To show that there are also non-Sidon sets sharing this property
we observe a simple lemma.

\begin{lemma}\label{lem9}
If $Z$ is a rich subspace of $X$, then $L_{1}(Z)$ is a rich subspace
of the Bochner space $L_{1}(X)$.
\end{lemma}

\begin{proof}
It is enough to check the definition of narrowness of the quotient
map on vector-valued step functions. Thus the assertion of the lemma
is reduced to the assertion that $Z \oplus_{1} \dots \oplus_{1} Z$
is a rich subspace of  $X \oplus_{1} \dots \oplus_{1} X$; but this
has been proved in \cite{BKSW}.
\end{proof}

Now if $\Lambda \subset \Z$ is a co-Sidon set, then
$L_{1}(L_{1,\Lambda}) \cong L_{1,\Z\times \Lambda} (\T^2)$ is a rich
subspace of $L_{1}(L_{1}) \cong L_{1}(\T^2)$, and $\Lambda'=\Z \times
(\Z\backslash \Lambda)$ is a non-Sidon set with reference to the
group $\T^2$ for which Corollary~\ref{cor6} is valid.

%--------------------------------------------------------
%%%%%%%%%%%%%%%%%%%%%%%%%%%%%%%%%%%%%%%%%%%%%%%%%%%%%%%%%
\section{Some examples of small but rich subspaces} 
\label{sec3}

In this section we provide examples of
nonseparable Banach spaces and separable rich subspaces. 

First we give a handy reformulation of richness. We let 
$$
D(x,y,\eps) = \{z\in X\dopu 
\|x+y+z\|> 2 - \eps,\  \|y+z\|< 1 + \eps \}
$$
for $x,y\in S(X)$. 

\begin{lemma}\label{66.1}
The following assertions are equivalent for a Banach space~$X$.
\begaeq
\item
$Z$ is a rich subspace of $X$.
\item
For every $x,y\in S(X)$ and every $\eps>0$,
$$
y\in \coq \bigl( y+
\bigl( D(x,y,\eps) \cap Z \bigr) \bigr).
$$
\item
For every $x,y\in S(X)$ and every $\eps>0$,
$$
0\in \coq 
\bigl( D(x,y,\eps) \cap Z \bigr) .
$$
\endaeq
\end{lemma}

\begin{proof}
(i) $\Leftrightarrow$ (ii) is a consequence of the Hahn-Banach theorem,
and (ii) $\Leftrightarrow$ (iii) is obvious.
\end{proof}

For $Z=X$, (ii) boils down to condition~(iii) of Lemma~\ref{l:MAIN}.

In the examples we are going to present $Z$ will be a space $C(K,E)$
embedded in a suitable space $X$. The type of space we have in mind
will be defined next.

\begin{definition}\label{66.2}
Let $E$ be a Banach space and $X$ be a sup-normed space  of bounded
$E$-valued functions on a  compact space~$K$. 
The space $X$ is said to be
a $C(K,E)$-superspace if it contains  $C(K,E)$ and for every 
$f \in X$, every $\eps>0$ and every open 
subset $U \subset K$
there exists an element $e \in E$, $\|e\| > (1 -\eps)\sup_U \|f(t)\|$, 
and a nonvoid open subset $V \subset U$ such that 
$\|e - f(\tau)\| < \eps$ for every $\tau \in V$.
\end{definition}

Basically, $X$ is a $C(K,E)$-superspace if every element of $X$ is
large and almost constant on suitable open sets.

Here are some examples of this notion.

\begin{prop}\label{66.3}
\mbox{}
\begsta
\item
$D[0,1]$, the space of bounded functions on $[0,1]$ that are
right-con\-tin\-uous and have left limits everywhere and are continuous at
$t=1$, is a $C[0,1]$-superspace.
\item
Let $K$ be a  compact Hausdorff space and $E$ be a Banach space.
Then $C_w(K,E)$, the space of weakly continuous functions
from $K$ into $E$, is a $C(K,E)$-superspace.
\endsta
\end{prop}

\begin{proof}
(a) $D[0,1]$ is the uniform closure of the span  of the  step functions
$\chi_{[a,b)}$, $0\le a<b<1$, and $\chi_{[a,1]}$, $0\le a<1$; hence
the result.

(b) Fix $f$, $U$ and $\eps$ as in Definition~\ref{66.2}; without loss
of generality we assume that $\sup_U \|f(t)\|=1$. 
Consider the open set $U_0 = \{t\in U\dopu \|f(t)\| > 1-\eps\}$. Now
$f(U_0)$ is  relatively weakly compact  since $f$ is weakly
continuous; hence it is dentable \cite[p.~110]{BL1}. Therefore there exists a
halfspace $H= \{x\in E\dopu x^*(x)>\alpha\}$ such that $f(U_0) \cap H$
is nonvoid and has diameter ${<\eps}$. Consequently, $V:= f^{-1}(H)
\cap U_0$ is an open subset of $U$ for which $\|f(\tau_1 ) -
f(\tau_2)\| < \eps$ for all $\tau_1, \tau_2 \in V$. This shows that
$C_w(K,E)$ is a $C(K,E)$-superspace. 
\end{proof}

The following theorem explains the relevance of these ideas.

\begin{thm}\label{66.4}
If $X$ is a $C(K,E)$-superspace and $K$ is perfect, then $C(K,E)$ is
rich in $X$; in particular, $X$ has the \DP.
\end{thm}

\begin{proof}
We wish to verify condition~(iii) of Lemma~\ref{66.1}. Let $f,g\in S(X)$
and $\eps>0$.
We first find an open set $V$  and an element $e\in E$,
$\|e\|>1-\eps/4$, such that $\|e-f(\tau)\| < \eps/4$ on~$V$. Given
$N\in\N$, find open nonvoid pairwise disjoint subsets $V_{1},\dots,
V_{N}$ of~$V$. Applying the definition again, we obtain elements
$e_{j}\in E$ and open subsets $W_{j}\subset V_{j}$ such that 
$\|e_{j}\| > (1-\eps/4) \sup_{V_{j}} \|g(t)\|$ and
$\|e_{j}-g(\tau)\|<\eps/4$ on $W_{j}$. Let $x_{j}= e-e_{j}$, let
$\varphi_{j}\in C(K)$ be a
positive function supported on $W_{j}$ of norm~$1$ and let $h_{j}=
\varphi_{j}\otimes x_{j}$. Now if $t_{j}\in W_{j}$ is selected to
satisfy $\varphi_{j}(t_{j})=1$, then
$$
\|f+g+h_{j}\| \ge \|(f+g+h_{j})(t_{j})\| > \|e+e_{j}+x_{j}\| -\eps/2
>2-\eps
$$
and 
$$
\|g+h_{j}\|<1+\eps
$$
since $\|g(t)+h_{j}(t)\|\le 1$ for $t\notin W_{j}$, and for $t\in W_{j}$
$$
\|g(t)+h_{j}(t)\|\le \|e_{j} + \varphi_{j}(t)x_{j}\| +\eps/4
\le (1-\varphi_{j}(t)) \|e_{j}\| + \varphi_{j}(t) \|e\| + \eps/4.
$$
This shows that $h_{j}\in D(f,g,\eps) \cap
C(K,E)$. But the supports of the $h_{j}$ are pairwise disjoint, hence 
$\|1/N \sum_{j=1}^N h_{j} \| \le 2/N \to 0$.
\end{proof}

\begin{cor}\label{66.5}
\mbox{}
\begsta
\item
$C[0,1]$ is a separable rich subspace of the nonseparable space
$D[0,1]$.
\item
If $K$ is perfect, then $C(K,E)$
is a rich subspace of $C_{w}(K,E)$. In particular, $C([0,1],\ell_{p})$
is a separable rich subspace of the nonseparable space
$C_{w}([0,1],\ell_{p})$ if $1<p<\infty$.
\endsta
\end{cor}

Let us remark that there exist nonseparable spaces
with the Daugavet property with only nonseparable rich
subspaces. Indeed, an $\ell_\infty$-sum of
uncountably many spaces with the \DP\ is  an example of this
phenomenon. 
To see this we need the result from \cite{BKSW} that whenever $T$ is a
narrow operator on $X_1 \oplus_\infty X_2$, then the restriction of
$T$ to $X_1$ is narrow too, and in particular it is not bounded from
below. Now let $X_i$, $i\in I$, be Banach spaces with the \DP\ and let
$X$ be their $\ell_\infty$-sum. If $Z$ is a rich subspace of $X$, then
by the result quoted  above there exist elements $x_i \in S(X_i)$ and
$z_i\in Z$ with $\|x_i - z_i\|\le 1/4$; hence $\|z_i - z_j\|\ge 1/2$
for $i\neq j$. If $I$ is uncountable, this implies that $Z$ is
nonseparable.

%--------------------------------------------------------
%%%%%%%%%%%%%%%%%%%%%%%%%%%%%%%%%%%%%%%%%%%%%%%%%%%%%%%%%
\section{The Daugavet property and tensor products}
\label{sec4}

One may consider the space  $C(K,E)$ as the injective
tensor product of  $C(K)$ and $E$; see for instance \cite[Ch.~VIII]{DiUh} 
or \cite[Ch.~3]{Ryan}
for these matters. It is known that $C(K,E)$ has the \DP\ whenever
$C(K)$ has, regardless of $E$ (\cite{Kadets} or \cite{KadSSW}), and
it is likewise true that $C(K,E)$ has the \DP\ whenever $E$ has,
regardless of $K$ \cite{MarPaya}. This raises
the natural question whether the injective tensor product of
two spaces has the Daugavet property if at least one 
factor has. 

We first give a positive answer for the class of rich subspaces of
$C(K)$; for example, a
uniform algebra is a rich subspace of $C(K)$ if $K$ denotes its Silov
boundary and is perfect.

\begin{prop}\label{prop8.1}
If $X$ is a rich subspace of some $C(K)$-space, then 
$X \widehat{\otimes}_{\eps} E$, the completed injective tensor product
of $X$ and $E$, is a rich subspace of $C(K) \widehat{\otimes}_{\eps}
E$  for every Banach space~$E$; in particular, it has the \DP. 
\end{prop}

\begin{proof}
We will consider $X \widehat{\otimes}_{\eps} E$ as a subspace of
$C(K,E)$. In order to verify (iii) of Lemma~\ref{66.1} let $f,g\in 
S( C(K,E) )$ and $\eps>0$ be given. Further,
let $\eta>0$ be given. We wish to construct functions
$h_1,\dots,h_n\in D(f,g,\eps) \cap X \widehat{\otimes}_{\eps} E$ 
such that $\|\frac1n \sum_{j=1}^n h_j\| \le 2\eta$. 

There is no loss in assuming that $\eta\le\eps$. Consider
$U=\{t\dopu \|f(t)\|>1-\eta/2\}$. By reducing $U$ if necessary we may
also assume that $\|g(t)-g(t')\|<\eta$ for $t,t'\in U$. Fix
$n\ge 2/\eta$  and pick $n$ pairwise disjoint open nonvoid subsets
$U_{1},\dots,U_{n}$ of~$U$; this is possible since $K$ must be
perfect, for $C(K)$ carries a narrow operator, viz.\ the quotient map
$q\dopu C(K)\to C(K)/X$. By applying 
\cite[Th.~3.7]{KadSW2} 
%or by Proposition~\ref{Cvs.nar},
to $q$ we infer that there exists,
for each~$j$, a function $\psi_{j}\in X$ 
with $\psi_{j}\ge0$, $\|\psi_{j}\|=1$
and $\psi_{j} <\eta/2$ off~$U_{j}$. 
Choose $t_{j}\in U_{j}$ with $\psi_{j}(t_{j})=1$. We define
$$
h_{j}=  \psi_{j}\otimes (f(t_{j}) - g(t_{j})) 
\in X \widehat{\otimes}_{\eps} E 
$$
and claim that $h_j\in D(f,g,\eta) \subset D(f,g,\eps)$.
In fact, 
$$
\|f+g+h_j\| \ge \|f(t_j) + g(t_j) + h_j(t_j) \| = 2\|f(t_j)\| >
2-\eta.
$$
Also, $\|g+h_j\| < 1+\eta$, for if $t\in U_{j}$, then 
\begin{eqnarray*}
\|g(t) + h_{j}(t)\| 
&\le&
\|g(t_j) + h_{j}(t)\| + \|g(t) -g(t_j)\|  \\
&<&
\|(1-\psi_{j}(t)) g(t_j) + \psi_{j}(t) f(t_{j})\| + \eta \\
&\le&
1+\eta,
\end{eqnarray*}
and for $t\notin U_{j}$  we clearly have $\|g(t) + h_{j}(t)\| <  1+\eta$.

It is left to estimate $\|\frac1n \sum_{j=1}^n h_j\|$.
If $t$ does not belong to any of the $U_{j}$, we have
$$
\biggl\|  \frac1n \sum_{j=1}^n h_{j}(t) \biggr\| \le \eta,
$$
and if $t\in U_{i}$, we have
$$
\biggl\|  \frac1n \sum_{j=1}^n h_{j}(t) \biggr\| 
\le \frac{n-1}n \eta + \frac1n \|h_i(t)\| 
\le \eta + \frac2n \le 2\eta
$$
by our choice of $n$.
\end{proof}

In general, however, the above question has a negative answer.

\begin{thm}\label{theo8.2}
There exists a two-dimensional complex Banach space $E$ such that
$L_{1}^{\C}[0,1] \widehat{\otimes}_{\eps} E$  fails the \DP, where
$L_{1}^{\C}[0,1]$ denotes the space of complex-valued $L_{1}$-functions.
\end{thm}

\begin{proof}
Consider the subspace $E$ of complex $\ell_\infty^{6}$ spanned by
the vectors $x_1=(1,1,1,1,1,0)$ and 
$x_2 = (0, \frac12, -\frac12, \frac{i}2,-\frac{i}2, 1)$.
The injective tensor product of  $E$ and $L_1^{\C}[0,1]$ can be identified 
with the space of 6-tuples of functions
$f=(f_1, \dots , f_6)$ of the form 
$g_1 \otimes x_1 + g_2 \otimes x_2$, $g_1, g_2 \in  L_1^{\C}[0,1]$, 
with the norm $\|f\|= \max_{k=1, \dots , 6}\|f_k\|_1$. 
To show that this space does not have the Daugavet property,
consider the slice 
$$
S_\eps = \biggl\{ f=(f_1, \dots , f_6) \in   L_1^{\C}[0,1]\otimes E\dopu  \re
\int_0^1 f_1(t)\, dt > 1 - \eps, \ \|f\| \le 1 \biggr\}.
$$
Every $f=g_1 \otimes x_1 + g_2 \otimes x_2 \in S_\eps$ satisfies the
conditions 
$$
\|g_1\| > 1 - \eps, \quad \max\{\|g_1 \pm \textstyle \frac12
g_2\|,  \|g_1 \pm \frac{i}2 g_2\|\} \le 1.
$$ 
Now the complex space $L_{1}$ is complex uniformly convex 
\cite{Glob}. Therefore,
there exists a function 
$\delta (\eps)$, which tends to 0 when $\eps$ tends to 0, such that 
$\|g_2\| < \delta (\eps)$ for every 
$f=g_1 \otimes x_1 + g_2 \otimes x_2 \in S_\eps$. This implies that
for every $f \in S_\eps$ 
$$
\|1 \otimes x_{2} + f\| \le \frac32 + \delta (\eps). 
$$
So if $\eps$ is small enough, there is no $f \in S_\eps$ 
with $ \|1 \otimes x_{2} + f\| > 2 - \eps$. 
By Lemma~\ref{l:MAIN}, this proves
that this injective tensor product does not have the Daugavet property.
\end{proof}

For the projective norm it is known that $L_1(\mu)
\widehat{\otimes}_{\pi } E = L_{1}(\mu,E)$  has the \DP\ regardless
of $E$ whenever $\mu$ has no atoms \cite{KadSSW}. Again, there is a
counterexample in the general case.

\begin{cor}\label{cor8.3}
There exists a two-dimensional complex Banach space $F$ such that
$L_{\infty}^{\C}[0,1] 
\widehat{\otimes}_{\pi} F$      fails the \DP, where
$L_{\infty}^{\C}[0,1]$ 
denotes the space of complex-valued $L_{\infty}$-functions.
\end{cor}

\begin{proof}
Let $E$ be the two-dimensional space from Theorem~\ref{theo8.2};
note that $( L_{1}^{\C} \widehat{\otimes}_{\eps} E)^* = L_{\infty}^{\C}
\widehat{\otimes}_{\pi} E^*$. Since the \DP\ passes from a dual space
to its predual, $F:=E^*$ is the desired example.
\end{proof}

%%%%%%%%%%%%%%%%%%%%%%%%%%%%%%%%%%%%%%%%%%%%%%%%%%%%%%%%%%%%%%%%%%%%%%
\section{Questions}
\label{secQ}

We finally mention two questions that were raised by
A.~Pe{\l}czy\'nski which we have not been able to solve.

(1) Is there a rich subspace of $L_1$ with the Schur property?
It was recently proved in \cite{KadWer} that the  subspace
 $X\subset L_1$ constructed by 
Bourgain and Rosenthal in \cite{BourRos}, which has the Schur property
and fails the RNP, is a space with the \DP; however, it is not rich in~$L_1$.

(2) If $X$ is a subspace of $L_{1}$ with the RNP, does $L_{1}/X$ have the
\DP? The answer is positive for reflexive spaces \cite{KadSSW}, for
$H^1$ \cite{Woj92} and a certain space constructed by Talagrand 
\cite{Tala5} in his (negative) 
solution of the three-space problem for $L_{1}$ \cite{KadSSW}.

%%%%%%%%%%%%%%%%%%%%%%%%%%%%%%%%%%%%%%%%%%%%%%%%%%%%%%%%%%%%%%%%%%%%%%
%\bibliography{az}
%\bibliographystyle{standard}

\end{document}